\documentclass[12pt]{amsart}
\usepackage[T2A]{fontenc}
\usepackage[cp1251]{inputenc}
\usepackage{amsmath,amsfonts}
\usepackage{amssymb}
\usepackage{amscd}
\usepackage[dvips]{graphicx}

\newcommand{\im}{\ensuremath{\mbox{Im}\,}}
\newcommand{\re}{\ensuremath{\mbox{Re}\,}}

\newcommand{\defin}[1]{\bf \noindent Definition #1.\rm\,\,}
\newcommand{\prop}[1]{\bf \noindent Proposition #1.\it\,\,}

\newcommand{\theor}[1]{\noindent \bf Theorem #1.\it\,\,}
\newcommand{\doc}{\noindent \bf Proof.\rm\,\,}
\newcommand{\remark}[1]{\noindent \bf Remark #1.\rm\,\,}
\newcommand{\CC}[1]{\mathbb{C}^{#1}}
\newcommand{\RR}[1]{\mathbb{R}^{#1}}

\newcommand{\dww}{\frac{\partial}{\partial w_2}}
\newcommand{\dwww}{\frac{\partial}{\partial w_3}}
\newcommand{\dz}{\frac{\partial}{\partial z}}
\newcommand{\dzz}{\frac{\partial}{\partial z_2}}
\newcommand{\dzzz}{\frac{\partial}{\partial z_3}}

\newcommand{\x}{X_1}
\newcommand{\xp}{X'_1}
\newcommand{\xx}{X_2}
\newcommand{\xxx}{X_3}
\newcommand{\xxxx}{X_4}\newcommand{\lr}{\longrightarrow}

\newcommand{\dimm}[1]{\mbox{dim}\,\mathfrak{aut}\,{#1}_p}
\newcommand{\dimp}[1]{\mbox{dim}\,\mathfrak{aut}_p{#1}_p}
\newcommand{\g}{\mathfrak{g}(M)}
\newcommand{\type}[1]{\medskip\noindent\bf\mbox{Type}\,{#1}.\rm}

\title{Classification of homogeneous CR-manifolds in dimension 4}
\author {V.K.Beloshapka and I.G.Kossovskiy}

\begin{document}
\maketitle

\begin{abstract}
Locally homogeneous CR-manifolds in dimension 3 were classified,
up to local CR-equivalence, by E.Cartan. We classify, up to local
CR-equivalence, all locally homogeneous CR-manifolds in dimension
4. The classification theorem enables us also to classify all
symmetric CR-manifolds in dimension 4, up to local biholomorphic
equivalence. We also prove that any 4-dimensional real Lie algebra
can be realized as an algebra of affine vector fields in a domain
in $\CC{3}$, linearly independent at each point.
\end{abstract}

\section{introduction}

Let $M$ be a real-analytic CR-manifold, generically embedded to
the complex space $\CC{n+k},\,n=\mbox{CR\,dim}\,M,
\,k=\mbox{codim}\,M.$ An important class of CR-manifolds is the
class of homogeneous CR-manifolds. A CR-manifold $M$ is called \it
homogeneous, \rm if its CR-automorphism group acts transitively on
it. The CR-manifold $M$ is called \it locally homogeneous, \rm  if
germs of $M$ at any two points are CR-equivalent. In other words,
locally homogeneous CR-manifolds are that ones which are "the same
at all points" up to local biholomorphic transformations (note
that, due to the classical result of Tomassini\,\cite{tomassini},
the local CR-equivalence is equivalent to the local biholomorphic
equivalence in the real-analytic category, so we don't make a
difference in these notions in what follows). A very useful
equivalent definition of local homogeneity is as follows (see
\cite{zaitsev} for possible equivalent definitions of local
homogeneity). A \it holomorphic vector field on $M$ at a point
$p\in M$ \rm is a vector field, which is tangent to $M$ at each
point and has the form
$$2\re\left(f_1(z)\frac{\partial}{\partial z_1}+...+f_{n+k}(z)\frac{\partial}{\partial z_{n+k}}\right),$$  where $f_j(z)$ are holomorphic functions in a neighborhood of $p$ in the ambient space. In what follows we skip the operator $2\re(\cdot)$. Holomorphic vector fields at $p$ are exactly that ones which generate a local flow of biholomorphic transformations at $p$, preserving $M$. These vector fields constitute a Lie algebra with respect to the Lie bracket of vector fields, which is called \it the infinitesimal automorphism algebra of $M$ at $p$ \rm and is denoted by $\mathfrak{aut}\,M_p$.  The \it stability subalgebra \rm of $M$ at the point $p$ is the subalgebra
$\mathfrak{aut}_p\,M_p\subset\mathfrak{aut}\,M_p$, which consists
of vector fields, vanishing at $p$. This algebra generates the \it
stability group \rm (or the \it isotropy group) \rm of $M$ at $p$,
which consists of biholomorphic automorphisms of the germ of $M$
at $p$, preserving the fixed point $p$. An \it evaluation mapping
\rm is the natural mapping
$\varepsilon_p:\,\mathfrak{aut}\,M_p\rightarrow T_p M$, given by
the formula $\varepsilon_p(X)=X|_p$. $M$ is called  \it locally
homogeneous at \rm $p$, if the mapping $\varepsilon_p$ is
surjective (i.e. the values of vector fields from
$\mathfrak{aut}\,M_p$ at $p$ form all the tangent space $T_p M$).
$M$ now is called \it locally homogeneous, \rm if it is locally
homogeneous at all points.

Locally homogeneous CR-manifolds in dimension 3 were classified by
E.Cartan\,\cite{cartan}. Any such manifold is locally
CR-equivalent to one of the following hypersurfaces in $\CC{2}$:
the hyperplane $\im z_2=0$ (the case $\dimp{M}=\infty$), the
hypersphere $|z_1|^2+|z_2|^2=1$ (the case $\dimp{M}=5$) or one of
the \it Cartan's homogeneous surfaces \rm (the case $\dimp{M}=0$).
Note that for a locally homogeneous CR-manifold we clearly have
$\dimm{M}=\dimp{M}+\mbox{dim}\,M$, so in the mentioned three cases
we have $\dimm{M}=\infty,\,\dimm{M}=8$ and $\dimm{M}=3$
correspondingly.  Cartan's classification is essentially based on
the Bianki's classification of real 3-dimensional Lie algebras.
Classification of locally homogeneous CR-manifolds in dimension 5
(which is essentially the classification of locally homogeneous
hypersurfaces in $\CC{3}$) is in progress. Partial results in the
Levi non-degenerate case have been obtained by A.Loboda, who
presented the desired classification for CR-manifolds with the
condition $\dimp{M}>1$ (see \cite{loboda1},\cite{loboda2}) and
partial classification in the case $\dimp{M}=1$ (this case is to
be completed soon; see \cite{loboda3} for details). Unfortunately,
in the case $\dimm{M}=0$ A.Loboda's approach is unapplicable and
this case is to be studied by different methods (see \cite{C3} for
some examples concerned with this case). In the Levi-degenerate
case the complete classification has been obtained by G.Fels and
W.Kaup in \cite{kaup} (see also \cite{kaup1}).

In the present paper we classify all locally homogeneous CR-manifolds in dimension 4, i.e. we do the next step after E.Cartans
classification. Since for a generic CR-submanifold $\mbox{dim}\,M=2n+k$ holds, from $k,n>0$ we get $n=1,k=2$
(the classification of locally homogeneous CR-manifolds in the case of a complex manifold ($k=0$) and in the case of a totally real manifold ($n=0$) is trivial: in the first case $M$ is locally CR-equivalent to $\CC{2}$, in the second case $M$ is locally
CR-equivalent to a totally real 4-plane in $\CC{4}$). So in what follows $M$ is supposed to be a real-analytic locally homogeneous
CR-manifold of dimension 4 and codimension 2, generically embedded to the complex space $\CC{3}$.

A natural example of homogeneous CR-manifolds in the case under consideration is as follows. Let $L$ be an affinely homogeneous curve in $\RR{3}$. Then the tube manifold $M=\{z\in\CC{3}:\,\,\im z\in L\}$ is a
homogeneous 4-dimensional CR-manifold. The homogeneity is provided by the abelian group of real translations $z\mapsto z+a,\,a\in\RR{3}$, and by the one-dimensional affine group $z\mapsto A(t)z+b(t)$, where $y\mapsto A(t)y+b(t)$ is the one-dimensional affine group, providing the homogeneity of $L$.

Another important example of a homogeneous 4-dimensional CR-manifold is concerned with the \it main trichotomy \rm for
CR-manifolds under consideration, demonstrated by V.Beloshapka, V.Ezhov and G.Schmalz in \cite{small},\cite{bes}. To formulate it, we firstly introduce the notion of total non-degeneracy.
A 4-dimensional CR-manifold $M$ in $\CC{3}$ is called \it totally non-degenerate at a point $p$, \rm if
$$T_p^{\CC{}}M+[T^{\CC{}}M,T^{\CC{}}M]_p+[T^{\CC{}}M,[T^{\CC{}}M,T^{\CC{}}M]]_p=T_pM,$$
where $T^{\CC{}}M$ is the bundle of complex tangent planes to $M$;
otherwise we call $M$ just \it degenerate \rm at $p$. At a generic
point this non-degeneracy condition is equivalent to the \it
holomorphic non-degeneracy \rm (see \cite{bouendy}). A totally
non-degenerate at each point CR-manifold $M$ is also called \it
Engel-type manifold \rm (see \cite{bes2}). As it was shown in
\cite{bes}, the total non-degeneracy condition is also equivalent
to the existence of local holomorphic coordinates
$(z,w_2,w_3)\in\CC{3}$, in which $p$ is the origin and $M$ is
given as
$$\im w_2=|z|^2+O(3),\,\im w_3=2|z|^2\re z+O(4),$$ where $z,w_2,w_3$ are assigned the weights
$$[z]=1,\,[w_2]=2,],[w_3]=3$$ and $O(3),O(4)$ are terms of weights $\geq 3$ and $\geq 4$ correspondingly. The manifold, for which
the remainders vanish, is called \it the 4-dimensional CR-cubic
\rm (we will call it just \it the cubic \rm and denote by $C$).
The cubic is the main example of a homogeneous totally
non-degenerate CR-manifold. The cubic is a particular case of a
\it model manifold \rm (see, for example, \cite{univers}) and has
many remarkable properties (see
\cite{small},\cite{bes},\cite{bes2},\cite{C3} for details), the
main one is given by the following trichotomy for a 4-dimensional
locally homogeneous CR-submanifold $M$ in $\CC{3}$:

(1) $\mbox{dim}\,\mathfrak{aut}_p M_p=\infty$, which occurs if and
only if $M$ is locally biholomorphically equivalent to a direct
product $M^3\times \RR{1}$, where $M^3\subset\CC{2}_{z_1,z_2}$ is
a locally homogeneous hypersurface in $\CC{2}$,
$\RR{1}\subset\CC{1}_{z_3}$ is a real line (the totally degenerate
case).

(2)  $\mbox{dim}\,\mathfrak{aut}_p M_p=1$, which occurs if and only if $M$ is locally biholomorphically
equivalent to the cubic $C$.

(3) $\mbox{dim}\,\mathfrak{aut}_p M_p=0$ for all other totally non-degenerate manifolds (the \it rigidity phenomenon). \rm

As it was demonstrated in \cite{bes2}, the cubic can be also
singled out among totally non-degenerate surfaces as the
(essentially) unique \it flat \rm surface with respect to some
special CR-curvature.

It is seen from the above trichotomy that \it the homogeneity of a 4-dimensional totally non-degenerate CR-manifold, which is locally non-equivalent to the cubic, is provided by a 4-dimensional algebra of holomorphic vector fields, namely by its infinitisemal automorphism algebra. \rm
It also follows from the trichotomy that \it a local  CR-mapping between two 4-dimensional non-degenerate CR-manifolds, locally
non-equivalent to the cubic, induces a local biholomorphic mappings between the vector field algebras, providing the homogeneity
of the manifolds \rm (see proposition 3.2 for details).
These two observation essentially simplifies the desired classification, which is finally given by the following theorem
(we use notations, associated to the cubic and denote the coordinates in $\CC{3}$ by $z,w_2,w_3$ and also set
$z=x+iy,\,w_j=u_j+iv_j$):

\smallskip

\bf Main Theorem. \it Any real-analytic  4-dimensional locally
homogeneous CR-submanifold $M$ in $\CC{3}$ is locally
CR-equivalent to one of the following pairwise locally
CR-inequivalent homogeneous surfaces:

\smallskip

CASE 1 - $\mbox{dim}\,\mathfrak{aut}_p M_p=\infty$ (the degenerate
case):

\smallskip

(1.1) $v_2=0,\,v_3=0$ (the real plane).

\smallskip

(1.2) $|z|^2+|w_2|^2=1,\,v_3=0$.

\smallskip

(1.3) $M^3\times \RR{1}$, where $M^3\subset\CC{2}_{z,w_2}$ is one
of the Cartan's homogeneous surfaces in $\CC{2}$,
$\RR{1}\subset\CC{1}_{w_3}$ is a real line.

\medskip

CASE 2 - $\mbox{dim}\,\mathfrak{aut}_p M_p=1$:

\smallskip

(2.1) $\{v_2=|z|^2,\,v_3=2|z|^2\re z\}\,\sim\{v_2=y^2,\,v_3=y^3\}$ (the cubic).

\medskip

CASE 3 - $\mbox{dim}\,\mathfrak{aut}_p M_p=0$:\rm

\smallskip

(3.1) $v_2=xe^y+\gamma
ye^y,\,v_3=e^{y},\,\gamma\in\RR{}.$

\smallskip

(3.2) $v_2=\frac{x}{y}+\gamma\ln
 y,\,v_3=\frac{1}{y},\,\gamma\in\RR{}.$

\smallskip

(3.3) $v_2=xy^{\alpha}+\gamma
y^{\alpha+1},\,v_3=y^{\alpha},\,|\alpha|>1,\,\alpha\neq
2,\gamma\in\RR{}.$

\smallskip

(3.4) $v_2=xy\ln y+\gamma y^2,\,v_3=y\ln y,\,\gamma\in\RR{}$.

\smallskip

(3.5) $v_2=x\sqrt{1-y^2}+\gamma\mbox{arcsin}\,y,\,v_3=\sqrt{1-y^2},\,\gamma\in\RR{}$.

\smallskip

(3.6) $v_2=xv_3+\gamma (v_3^2+y^2),\,\exp
\left(q\,\mbox{arctg}\frac{v_3}{y}\right)=v_3^2+y^2,\,q>0,\gamma\in\RR{}$.

\smallskip

(3.7) $v_2=e^y,\,v_3=e^{x+\delta y},\,\delta\in\RR{}$.

\smallskip

(3.8) $v_2=e^{x+\alpha y}\cos \beta y,\,v_3=e^{x+\alpha y}\sin
\beta y,\,\beta > 0,(\alpha,\beta)\neq(0,1).$

\smallskip

(3.9) $v_2=e^x y \cos y,\,v_3=e^{x} y \sin y.$

\smallskip

(3.10)
$v_2=y^{\alpha},\,v_3=y^{\beta},\,1<\alpha<\beta,\,(\alpha,\beta)\neq(2,3).$

\smallskip

(3.11) $v_2=e^{ay},\,v_3=e^{y},\,0<|a|<1.$

\smallskip

(3.12) $v_2=\mbox{\rm ch}\,y,\,v_3=\mbox{\rm sh}\,y.$

\smallskip

(3.13) $v_2=y\ln y,\,v_3=y^{\alpha},\,\alpha\neq \{0;1\}.$

\smallskip

(3.14) $v_2=ye^{y},\,v_3=e^{y}.$

\smallskip

(3.15) $v_2=y^{2},\,v_3=e^{y}.$

\smallskip

(3.16) $v_2=y\ln^{2}y,\,v_3=y\ln y.$

\smallskip

(3.17) $v_2=e^{y}\cos\beta y,\,v_3=e^{y}\sin\beta y,\,\beta>0.$

\smallskip

(3.18) $v_2=y^{\alpha}\cos(\beta \ln y),\,v_3=y^{\alpha}\sin(\beta
\ln y),\,\beta>0.$

\smallskip

(3.19) $v_2=\cos y,\,v_3=\sin y$.

\smallskip

\smallskip

\rm

As a corollary we get the following statement.  \smallskip

\theor{1.2} Any locally homogeneous totally non-degenerate
CR-manifold is locally CR-equivalent to an affinely homogeneous
one. \rm\smallskip

Note that the last property holds for W.Kaup's list of
2-nondegenerate hypersurfaces in $\CC{3}$ but does not hold
E.Cartans and A.Loboda's lists.\smallskip

Theorem 1.2 implies the following "representation-type"\,
result.\smallskip

\theor{1.2'} Any 4-dimensional real Lie algebra can be realized as
an algebra of affine vector fields in a domain in $\CC{3}$,
linearly independent at each point.\rm\smallskip

The Main Theorem also allows us to formulate the following classification theorem.\smallskip

\theor{1.3} Any locally homogeneous totally non-degenerate
CR-manifold with non-trivial stability group is locally
CR-equivalent to one of the following pairwise locally
CR-inequivalent homogeneous surfaces: (2.1) (the cubic), (3.2),
(3.5), (3.12), (3.19).

In the first case the stability group at the origin looks as
$$z\mapsto\lambda z,\,w_2\mapsto\lambda^2 w_2,\,w_3\mapsto\lambda^3 w_3,\,\lambda\in\RR{*}$$ and thus is isomorphic to
$\RR{*}$; in all other cases the stability group is of $\mathbb{Z}_2$ - type and is generated by the authomorphism $$z\mapsto -z,\,w_2\mapsto -w_2,\,w_3\mapsto w_3$$ at the point $(0,0,i)$ for the surfaces (3.5), by the authomorphism $$z\mapsto -z,\,w_2\mapsto w_2,\,w_3\mapsto -w_3$$ at the origin for the surfaces (3.12), (3.19) and by the authomorphism $$z\mapsto w_3,\,w_2\mapsto -w_2+zw_3+1,\,w_3\mapsto z$$ at the point $(i,0,i)$ for the surfaces (3.2).\rm

\smallskip

Let $M$ be a Riemannian CR-manifold. $M$ is called \it Hermitian
CR-manifold, \rm   if the  Riemannian metric is compatible with
the almost complex structure $J$ on $M$ in the sense that
$||Jv||_p=||v||_p$ holds for each $v\in T^{\CC{}}_pM,p\in M$. A
Hermitian CR-manifold $M$ is called \it CR-symmetric, \rm  if for
each point $p\in M$ there exists a CR-isometry $s_p$ of $M$,
preserving the point $p$ and such that the differential $ds_p$,
restricted on the subspace $T^{\CC{}}_pM\oplus T^{TR}_p M\subset
T_pM,$  is minus identical. Here $T^{TR}_pM$ is the \it totally
real part of $T_pM$, \rm i.e. the orthogonal complement to the
subspace, spanned by $T^{\CC{}}_pM$ and by the values at $p$ of
arbitrary order Lie brackets of vector fields $X\in TM$ with the
condition $X_a\in T^{\CC{}}_aM$ for each $a\in M$ (for finite type
CR-manifolds, in particular for totally non-degenerate manifolds,
the subspace $T^{TR}_pM$ is trivial; for infinite type
CR-manifolds it is non-trivial and the condition for $ds_p$ to be
minus identical on $T^{TR}_pM$ guarantees the uniqueness of the
involution $s_p$). In the paper \cite{kaupzaitsev} some beautiful
connections between symmetric CR-manifolds and Hermitian symmetric
spaces are demonstrated (see also \cite{nachin}). In particular,
it is proved that any symmetric CR-manifold is also
CR-homogeneous.  The cubic $C$ is given in \cite{kaupzaitsev} as
an example of a symmetric CR-manifold in dimension 4. Using
theorem 3.1  and the Cartan's classification theorem, we obtain
the classification of \it all \rm symmetric CR-manifolds in
dimension 4.

\smallskip

\theor{1.4} Any symmetric CR-manifold of dimension 4 is locally
CR-equivalent to one of the following pairwise locally
CR-inequivalent symmetric surfaces:

\smallskip

CASE 1 - Degenerate Levi-flat manifolds:

 \smallskip

(1.a) $v_2=0,\,\,v_3=0$.

\smallskip

CASE 2 - Degenerate non Levi-flat manifolds:

\smallskip

(2.a) $v_2=y^2,\,\,v_3=0$.

\smallskip

(2.b) $y^2+v_2^2=1,\,\,v_3=0$.

\smallskip

(2.c) $y^2-v_2^2=1,\,\,v_3=0$.

\smallskip

(2.d) $1+|z|^{2}+|w_2|^{2}=a|1+z^{2}+w_2^{2}|,\,\,v_3=0,\,\,a>1.$

\smallskip

(2.e) $1+|z|^{2}-|w_2|^{2}=a|1+z^{2}-w_2^{2}|,\,\,v_3=0,\,\,a>1.$

\smallskip

(2.f) $-1+|z|^{2}+|w_2|^{2}=a|-1+z^{2}+w_2^{2}|,\,\,v_3=0,\,\,0<|a|<1.$

\smallskip

CASE 3 - Totally non-degenerate manifolds:

\smallskip

(3.a) $v_2=y^2,\,\,v_3=y^3$.

\smallskip

(3.b) $v_2=\cos y,\,\,v_3=\sin y$.

\smallskip

(3.c) $v_2=\mbox{\rm ch}\,y,\,\,v_3=\mbox{\rm sh}\,y$.

\smallskip

(3.d) $v_2=x\sqrt{1-y^2}+\gamma\mbox{\rm arcsin}\,y,\,y^2+v_3^2=1,\,\gamma\in\RR{}$.

\smallskip

(3.e) $v_2=\frac{x}{y}+\gamma\ln y,\,\,v_3=\frac{1}{y},\,\,\gamma\in\RR{}$.\rm

\smallskip

The paper is organized as follows. In section 2, using the main trichotomy, we consider locally homogeneous 4-dimensional CR-manifolds in $\CC{3}$ as orbits of the natural action of 4-dimensional Lie algebras of holomorphic vector fields in $\CC{3}$, and partially classify the orbits under the assumption that they are totally non-degenerate. In section 3 we analyze the homogeneous
4-dimensional CR-submanifolds in $\CC{3}$, obtained in section 2, and study the local CR-equivalence relations among them,
using the machinery of normal forms (see \cite{chern},\cite{bes2}) and the main trichotomy. This finally allows us to give a complete list of totally non-degenerate homogeneous manifolds under consideration and thus to prove the Main Theorem and then the theorems 1.3 and 1.4.

\remark{1.5} It follows from theorems 1.3 and 1.4 that the holomorphic authomorphism groups of the homogeneous surfaces (2.1) - (3.19) coincide with $\exp(\g)$ except the cases of symmetric surfaces, when the authomorphism group is a semidirect product
of $\exp(\g)$ and the stability group of a fixed point, described in theorem 1.3. Since in all cases $\g$ can be easily integrated,
this gives a description of the holomorphic authomorphism groups of the homogeneous surfaces (2.1) - (3.19).

\remark{1.6} Note that, in the same way as E.Cartans, A.Loboda's,
G.Fels and W.Kaups's lists of homogeneous CR-manifolds, our list
consists of \it globally \rm homogeneous surfaces, which equations
are given by elementary functions. Also note that, in the same way
as E.Cartans and G.Fels and W.Kaups's lists, (the non-degenerate
part of) our list consists of one "model" object with
positive-dimensional stability subalgebra and "rigid" objects with
trivial stability subalgebra (the "rigidity phenomenon").

\remark{1.7} The consideration of just real-analytic generically
embedded CR-submanifolds in $\CC{N}$ in this paper is motivated by
the fact that a CR-manifold, admitting a local transitive action
of a Lie group, is automatically real-analytic, and by the fact
that any real-analytic CR-manifold can be locally generically
embedded to $\CC{N}$ (see \cite{zaitsev}). This makes our
classification results general enough to be reformulated for
arbitrary (abstract) CR-manifolds.

The authors would like to thank Mike Eastwood  and Alex Isaev from the Australian National University
for useful discussions.

\section{Homogeneous CR-manifolds and 4-dimensional Lie algebras of holomorphic vector fields in $\CC{3}$}

In what follows M is supposed to be a real-analytic 4-dimensional
totally non-degenerate locally homogeneous generic CR-submanifold
in $\CC{3}$, $p=(p_1,p_2,p_3)$ is a fixed point on $M$. As it follows from the
above discussion, the homogeneity of $M$ is provided by some
4-dimensional real Lie algebra of holomorphic vector fields, which
coincides with $\mathfrak{aut}\,M_p$ if $M$ is locally
non-equivalent to the cubic $C$, or is a subalgebra of the
5-dimensional Lie algebra $\mathfrak{aut}\,M_p$ in case when $M$
is locally CR-equivalent to $C$ (see \cite{C3} for precise description of $\mathfrak{aut}\,C$).
We denote this algebra by $\g$. $M$ is the orbit of the natural local action of $\g$ in $\CC{3}$ at the point
$p$. Also we denote by $\x,\xx,\xxx,\xxxx$ a basis of the Lie
algebra $\g$.


As the values of  $\x,\xx,\xxx,\xxxx$ at $p$ span the 4-dimensional linear space $T_pM$, these values are linearly independent
over $\RR{}$. The next proposition shows that the total non-degeneracy property gives stronger restriction on these values.

\prop{2.1} If the vector fields $\x,\xx,\xxx$ span a 3-dimensional
Lie subalgebra $\mathfrak{a}$ of $\g$, then their values at $p$
are linearly independent over $\CC{}$.

\doc Suppose that
$\mbox{rk}_{\CC{}}\{\x|_p\,,\xx|_p\,,\xxx\,|_p\}<3$. Let
$\mathfrak{a}^{\CC{}}$ be the complexification of $\mathfrak{a}$.
Consider the real action of $\mathfrak{a}$ in $\CC{3}$ as well as
the complex action of $\mathfrak{a}^{\CC{}}$ in $\CC{3}$ and
denote the orbits of these local actions  at $p$ by $N$ and $L$
correspondingly. Since the values of $\x,\xx,\xxx$ at $p$ are
linearly independent over $\RR{}$, $N$ is a real 3-manifold. The
inequality $\mbox{rk}_{\CC{}}\{\x|_p\,,\xx|_p\,,\xxx\,|_p\}<2$ is
impossible because $M$ is generic (and hence the complexification
of $T_pM$ must coincide with $\CC{3}$). So we have
$\mbox{rk}_{\CC{}}\{\x|_p\,,\xx|_p\,,\xxx\,|_p\}=2$, $L$ is a
2-dimensional complex manifold and $N\subset L$ is a real
submanifold. Consider $T^{\CC{}}_pN\subset T^{\CC{}}_pL$. Since
$N\subset M$ and $T^{\CC{}}_pM$ is of dimension 1, we conclude
that $T^{\CC{}}_pN=T^{\CC{}}_pM$, consequently
$T^{\CC{}}_pM\subset T^{\CC{}}_pL$ and the same holds for all
neighbor points of $N$, which is a contradiction with the
total non-degeneracy condition. Hence $\mbox{rk}_{\CC{}}\{\x|_p\,,\xx|_p\,,\xxx\,|_p\}=3$, as required. \qed\\

In this section we obtain a partial classification of the class of
CR-manifolds under consideration, considering them as orbits of
the natural local action of 4-dimensional real Lie algebras in
$\CC{3}$. We use the classification of 4-dimensional real Lie
algebras, given in \cite{kruch} (the results of \cite{kruch} are
also claimed, for example, in \cite{nester}). There are 22 types
of such algebras: 10 solvable decomposable ones, 10 solvable
indecomposable ones, and 2 non-solvable decomposable ones. Some
types contain  real parameters. For our purposes it will be more
convenient to single out five types of solvable algebras, which do
not contain a 3-dimensional abelian ideal (according to
\cite{kruch}, these are types $A_{4.8},A_{4.7},A_{4.9},
A_{2.2}\oplus A_{2.2}$ and $A_{4.10}$ correspondingly). We also
denote by type VI  all solvable algebras, which contain a
3-dimensional abelian ideal (according to \cite{kruch}, these are
algebras of types $A_{3.1}\oplus A_{1}, A_{2.2}\oplus A_{2.1},
A_{3.j}\oplus A_1,\,j=3,...9$ and $A_{4.j},\,j=1,...,6$) and
denote by types VII and VIII correspondingly the two non-solvable
algebras $\mathfrak{so}_{2,1}(\RR{})\oplus\RR{1}$ and
$\mathfrak{so}_{3}(\RR{})\oplus\RR{1}$. According to the
classification, we consider 8 cases depending on the type of the
Lie algebra $\g$.

\type{I} Lie algebras of this type have the following commuting
relations:
\begin{gather*} \notag [\x,\xx]=0,\,[\x,\xxx]=0,\,[\xx,\xxx]=\x,\\
[\x,\xxxx]=(q+1)\x,\,[\xx,\xxxx]=\xx,\,[\xxx,\xxxx]=q\xxx,\,|q|\leq 1.\end{gather*}

 Applying now for an algebra of type I proposition
2.1, we conclude, that the values of $\x,\xx,\xxx$ at $p$ are
linearly independent over $\CC{}$ and we can rectify the commuting
vector fields $\x,\xx$ simultaneously in some neighborhood of $p$,
so in this neighborhood we have: $\x=\dww,\,\xx=\dwww$ (the
notations are taken from the introduction). From the commuting
relations we then have $\xxx=a(z)\dz+(b(z)+w_3)\dww+c(z)\dwww$ for
some analytic functions $a(z),b(z),c(z),\,a(p_1)\neq 0$. Then we
firstly rectify the non-zero vector field $a(z)\dz$ and after that
make a variable change  of kind $w_3\lr w_3+C(z)$. Then in the new
coordinates
$$\x=\dww,\,\xx=\dwww,\,\xxx=\dz+(\tilde{b}(z)+w_3)\dww+(\tilde{c}(z)+C_z)\dwww.$$ Now taking $C(z)$ from the equation $C_z+\tilde{c}(z)=0$ we get
$\xxx= \dz+(\tilde{b}(z)+w_3)\dww,$ and after a transformation
$w_2\lr w_2+B(z)$ for a function $B(z)$, satisfying
$B_z+\tilde{b}(z)=0$, we finally get $\xxx=\dz+w_3\dww$.

Now from the commuting relations for $\xxxx$ it is not difficult
to verify that $\xxxx$ must have the form
$\xxxx=(qz+l)\dz+((q+1)w_2+mz+n)\dww+(w_3+m)\dwww,\,l,m,n\in\CC{}$.
Also from the fact that the vector fields $\x,\xx,\xxx$ are
tangent to $M$, we can conclude that $M$ is given by equations
\begin{gather}v_2=x\psi(y)+\tau(y),\,v_3=\psi(y)\end{gather} for some real-analytic
functions $\psi(y),\tau(y)$. To see that, we present $M$ in the
form $v_2=F(x,y,u_2,u_3),\,v_3=G(x,y,u_2,u_3)$, which is possible
since $i\x,i\xx$ are transversal to $M$, and get from
the tangency conditions
$F_{u_2}=F_{u_3}=G_{u_2}=G_{u_3}=G_x=0,F_x=G$. Since $\xxxx$ is
tangent to $M$, we get the following conditions for $\psi,\tau$:
\begin{gather*}(q+1)x\psi(y)+(q+1)\tau(y)+m_1y+m_2x+n_2=\\
qx\psi+qxy\psi_y+l_2x\psi_y+l_1\psi+l_2\tau_y+qy\tau_y;\\
\psi+m_2=l_2\psi_y+qy\psi_y.\end{gather*}

Here $l=l_1+il_2,m=m_1+im_2,n=n_1+in_2$. Note that the linear part
of $\psi$ and $\tau$ in (1) can be annihilated by a linear
transformation. Also note that the total non-degeneracy condition
requires $\psi(y)\neq ay+b$. Then, solving the equations on
$\psi,\tau$ and making a linear variable change with respect to
$z$, we can put:

Case $q=0$: $\psi=Ae^y,\,\tau=Bye^y,\,A\neq 0$. After a scaling we
may suppose $A=1$.

Case $q=-1$: $\psi=\frac{A}{y},\,\tau=B\ln y,\,A\neq 0$. After a
scaling we may suppose $A=1$.

Case $q\neq\{0,\pm 1\}$:
$\psi(y)=Ay^{\alpha},\,\tau(y)=By^{\alpha+1},\,A\neq 0$, where
$\alpha=\frac{1}{q}$. After a scaling we may suppose $A=1$.

It is interesting that the case $q=1$ can't occur for a totally
non-degenerate manifold, i.e. all orbits of the Lie algebra of
type I with $q=1$ in $\CC{3}$ are degenerate.

\prop{2.2} Any totally non-degenerate CR-manifold with $\g$ of
type I is locally CR-equivalent to one of the following surfaces:

\begin{gather*}
 \mbox{ Type I\,a:}\,\,  v_2=xe^y+\gamma
ye^y,\,v_3=e^{y},\,\gamma\in\RR{}.\\
 \mbox{ Type I\,b:}\,\,  v_2=\frac{x}{y}+\gamma\ln
 y,\,v_3=\frac{1}{y},\,\gamma\in\RR{}.\\
 \mbox{ Type I\,c:}\,\,  v_2=xy^{\alpha}+\gamma
y^{\alpha+1},\,v_3=y^{\alpha},\,|\alpha|>1,\,\gamma\in\RR{}.\end{gather*} \rm

The restriction on $\alpha$ follows from the condition
$|q|\leq 1$.

\type{II} Lie algebras of this type have the following commuting
relations:
\begin{gather*} \notag [\x,\xx]=0,\,[\x,\xxx]=0,\,[\xx,\xxx]=\x,\\
[\x,\xxxx]=2\x,\,[\xx,\xxxx]=\xx,\,[\xxx,\xxxx]=\xx+\xxx.\end{gather*}

Since the commuting relations among $\x,\xx,\xxx$ are the same as
in case I, we conclude that in appropriate coordinates these
vector fields have the form
$\x=\dww,\,\xx=\dwww,\,\xxx=\dz+w_3\dww.$ From the commuting
relations for $\xxxx$ it is also not difficult to verify that
$\xxxx=(z+l)\dz+(2w_2+\frac{1}{2}z^2+mz+n)\dww+(z+w_3+m)\dwww,\,l,m,n\in\CC{}$.
After a translation along $z$ (which do not change $\x,\xx,\xxx$)
we may suppose $l=0$. Since $\x,\xx,\xxx$ are the same as in case
I, we can also conclude that $M$ has the form (1). The tangency
conditions for the vector field $\xxxx$ have the form:

\begin{gather*}2x\psi+2\tau+xy+m_1y+m_2x=x\psi+xy\psi_y+y\tau_y;\\
\psi+y+m_2=y\psi_y.\end{gather*} Annihilating the linear part of
$\psi$ and $\tau$, we can put: $\psi(y)=Ay\ln
|y|,\,\tau(y)=By^2,\,A\neq 0$. After a scaling we may suppose
$A=1$. Thus we have proved the following proposition.

\prop{2.3} Any totally non-degenerate CR-manifold with $\g$ of
type II is locally CR-equivalent to one of the following surfaces:
$$\mbox{Type\,II:}\quad v_2=xy\ln y+\gamma
y^2,\,v_3=y\ln y,\,\gamma\in\RR{}.$$

\type{III} Lie algebras of this type have the following commuting
relations:
\begin{gather*} \notag [\x,\xx]=0,\,[\x,\xxx]=0,\,[\xx,\xxx]=\x,\\
[\x,\xxxx]=2q\x,\,[\xx,\xxxx]=q\xx-\xxx,\,[\xxx,\xxxx]=\xx+q\xxx,\,q\geq 0.\end{gather*}

Since the commuting relations among $\x,\xx,\xxx$ are the same as
in case I, we conclude that in appropriate coordinates these
vector fields have the form
$\x=\dww,\,\xx=\dwww,\,\xxx=\dz+w_3\dww.$ From the commuting
relations for $\xxxx$ it is also not difficult to verify that
$\xxxx=(qz-w_3+l)\dz+(2qw_2+\frac{1}{2}z^2-\frac{1}{2}w_3^2+mz+n)\dww+(z+qw_3+m)\dwww,\,l,m,n\in\CC{}$.
After a translation along $z$ (which do not change $\x,\xx,\xxx$)
we may suppose $m=0$. Now we consider two cases.

Case 1: $q=0$. Since $\x,\xx,\xxx$ are the same as in case
I, we can also conclude that $M$ has the form (1). The tangency
conditions for the vector field $\xxxx$ have the form:

\begin{gather*}
xy-u_3\psi+n_2=(-u_3+l_1)\psi+x\psi_y(-\psi+l_2)+\tau_y(-\psi+l_2);\\
y=\psi_y(-\psi+l_2).\end{gather*} Annihilating the linear part of
$\psi$ and $\tau$, we can put:
$\psi(y)=\sqrt{R^2-y^2},\,\tau(y)=B\,\mbox{arcsin}\frac{y}{R},\,R>0$.
After a scaling we may suppose $R=1$, and $M$ is finally given as
 $$v_2=x\sqrt{1-y^2}+\gamma \arcsin
y,\,v_3=\sqrt{1-y^2},\,\gamma\in\RR{}.$$

Case 2: $q>0$. In that case after a translation along $w_2$ we may
suppose $n=0$. Also we make a variable change, which linearizes
our vector field algebra: \begin{gather}z\mapsto z,\,w_2\mapsto
2w_2-zw_3,\,w_3\mapsto w_3.\end{gather} As a result we have
$\x=2\dww,\,
\xx=\dwww-z\dww,\,\xxx=\dz+w_3\dww,\,\xxxx=(qz-w_3+l)\dz+(2qw_2-lw_3)\dww+(qw_3+z)\dwww$.
Replacing $\xxxx$ by $\xxxx-l_1\xxx$, we may suppose that $\re
l=0$. Now after the translations $z\mapsto z+qib,\,w_3\mapsto
w_3-ib$, where $b=\frac{l_1}{q^2+1}$, we get
$\x=2\dww,\,\xx=\dwww-(z-qib)\dww,\,\xxx=\dz+(w_3+ib)\dww,\,
\xxxx=(qz-w_3)\dz+(2qw_2-lw_3)\dww+(z+qw_3)\dwww$ (we replace
$\xxxx$ by $\xxxx-bl_2\x$). Presenting now $M$ in the form
$v_2=F(x,y,u_2,u_3),\,v_3=G(x,y,u_2,u_3)$, from the tangency
conditions for $\x,\xx,\xxx$ we get
$F_{u_2}=G_x=G_{u_2}=G_{u_3}=0,\,F_x=G+b,\,F_{u_3}=-y+bq$. Then
$M$ is given as
$$v_2=x\psi(y)+\tau(y)-yu_3+bqu_3+bx,\,v_3=\psi(y).$$

The tangency conditions for $\xxxx$ imply
\begin{gather*}
2q(x\psi+\tau-yu_3+bqu_3+bx)-bu_3(q^2+1))=\\=(qx-u_3)\psi+
x\psi_y(qy-\psi)+
\tau_y(qy-\psi)-(qy-\psi)u_3-\\-y(qu_3+x)+bq(qu_3+x)+b(qx-u_3);\\
q\psi+y=\psi_y(qy-\psi).\end{gather*} From the second equation we
get $$\psi_y=\frac{q\psi+y}{qy-\psi}.$$ This is a first order
homogeneous equation. The general solution is $\exp
(2q\mbox{arctg}\frac{\psi}{y})=c^2(\psi^2+y^2)$. Hence we get
$\psi_c(y)=\frac{1}{c}\psi_1(cy)$ and after a scaling we may
suppose for the original manifold $M$ that $\psi(y)=\psi_1(y)$.

For $\tau$ from the first tangency condition we get
\begin{gather}\frac{\tau_y}{\tau}=\frac{2q}{qy-\psi(y)}.\end{gather} Hence the
general solution is $\tau_c(y)=c\tau_1(y)$. It is straightforward
to check that $\tau=c(\psi^2+y^2)$ is actually the general
solution of (3). Annihilating the pluriharmonic terms in the right
hand side of the defining equations of $M$ and replacing $-yu_3$
by $xv_3$ (since the difference is pluriharmonic), we can present
$M$ after a scaling in the following way:

$$v_2=xv_3+\gamma (v_3^2+y^2),\,\exp \left(2q\,\mbox{arctg}\frac{v_3}{y}\right)=v_3^2+y^2,\,q>0,\gamma\in\RR{}.$$

\prop{2.4} Any totally non-degenerate CR-manifold with $\g$ of
type III is locally CR-equivalent to one of the following
surfaces:
\begin{gather*} \mbox{Type\,IIIa:}\quad v_2=x\sqrt{1-y^2}+\gamma
\arcsin
y,\,v_3=\sqrt{1-y^2},\,\gamma\in\RR{}.\\
\mbox{Type\,IIIb:}\quad v_2=xv_3+\gamma (v_3^2+y^2),\,\exp
\left(2q\,\mbox{\rm
arctg}\frac{v_3}{y}\right)=v_3^2+y^2,\,q>0,\gamma\in\RR{}.\end{gather*}

\type{IV} Lie algebras of this type have the following commuting
relations:
\begin{gather*} \notag [\x,\xx]=0,\,[\x,\xxx]=0,\,[\xx,\xxx]=\xx,\\
[\x,\xxxx]=\x,\,[\xx,\xxxx]=0,\,[\xxx,\xxxx]=0.\end{gather*}

Using proposition 2.1, we can rectify the commuting vector fields
$\x,\xx$, so $\x=\dww,\,\xx=\dwww$. From the commuting relations
$$\xxx=a(z)\dz+b(z)\dww+(w_3+d(z))\dwww,\,a(z)\neq 0.$$ Now we
firstly rectify $a(z)\dz$, then make a variable change $w_2\lr
w_2+B(z),\,w_3\lr w_3+D(z)$ for functions $B(z),D(z)$, satisfying
$B_z+b(z)=0,\,D_z-D(z)+d(z)=0$ and finally get
$\xxx=\dz+w_3\dwww$.

Now from the commuting relations for $\xxxx$ it is not difficult
to verify that $\xxxx$ must have the form
$\xxxx=l\dz+(w_2+m)\dww+ne^z\dwww,\,l,m,n\in\CC{}$. Also from the
fact that the vector fields $\x,\xx,\xxx$ are tangent to $M$, we
can conclude that $M$ is given by equations
\begin{gather}v_2=\psi(y),\,v_3=e^x\tau(y)\end{gather} for some real-analytic
functions $\psi(y),\tau(y)$. To see that, we present $M$ in the
form $v_2=F(x,y,u_2,u_3),\,v_3=G(x,y,u_2,u_3)$ and get from
the tangency conditions
$F_{u_2}=F_{u_3}=G_{u_2}=G_{u_3}=F_x=0,G_x=G$. Since $\xxxx$ is
tangent to $M$, we get the following conditions on $\psi,\tau$:
\begin{gather*}\psi+m_2=l_2\psi_y;\\
n_1e^x\sin y+n_2e^x\cos y=l_1e^x\tau+l_2e^x\tau_y.\end{gather*}

Note that linear terms in $\psi$ and terms of kind $a\sin y+b\cos
y$ in $\tau$ can be annihilated by transformations of kind $w_2\lr
w_2+A(z),\,w_3\lr w_3+B(z)$. Also from the total non-degeneracy
$\psi,\tau\neq 0$. Then after scalings $\psi=e^y,\tau=e^{\delta
y}$. Thus we have proved the following proposition.

\prop{2.5} Any totally non-degenerate CR-manifold with $\g$ of
type IV is locally CR-equivalent to one of the following surfaces:
$$\mbox{Type\,IV:}\quad v_2=e^y,\,v_3=e^{x+\delta y},\,\delta\in\RR{}.$$

\type{V} Lie algebras of this type have the following commuting
relations:
\begin{gather*} \notag [\x,\xx]=0,\,[\x,\xxx]=\x,\,[\xx,\xxx]=\xx,\\
[\x,\xxxx]=\xx,\,[\xx,\xxxx]=-\x,\,[\xxx,\xxxx]=0.\end{gather*}

Having $\g$ of type V, we firstly rectify the commuting vector
fields $\x,\xx$ (using proposition 2.1), so $\x=\dww,\,\xx=\dwww$.
Then from the commuting relations
$\xxx=a(z)\dz+(w_2+b(z))\dww+(w_3+d(z))\dwww,\,a(z)\neq 0$. After
a rectification of $a(z)\dz$ and a variable change $w_2\lr
w_2+B(z),\,w_3\lr w_3+D(z)$ for functions $B(z),D(z)$, satisfying
$B_z+b(z)=0;\,D_z+d(z)=0$, we finally get
$\xxx=\dz+w_2\dww+w_3\dwww$. From the commuting relations for the
vector field $\xxxx$ we then get
$$\xxxx=l\dz+(me^z-w_3)\dww+(ne^z+w_2)\dwww,\,l,m,n\in\CC{}.$$ Presenting $M$
in the form $v_2=F(x,y,u_2,u_3),\,v_3=G(x,y,u_2,u_3)$, we get from
the tangency conditions for $\x,\xx,\xxx$:
$F_{u_2}=F_{u_3}=G_{u_2}=G_{u_3}=0,F_x=F,G_x=G$, which means that
$M$ is given as \begin{gather}
v_2=e^x\psi(y),\,v_3=e^x\tau(y).\end{gather} The tangency
conditions for $\xxxx$ now give \begin{gather*}
m_1e^x\sin y+m_2e^x\cos y-e^x\tau=l_1e^x\psi+l_2e^x\psi_y;\\
n_1e^x\sin y+n_2e^x\cos
y+e^x\psi=l_1e^x\tau+l_2e^x\tau_y.\end{gather*} Solving these
equations and annihilating trigonometric terms $a\sin y+b\cos y$
by transformations of kind $w_j\lr w_j+C_je^z$ in (5), we get:

Case $l=l_1+il_2\neq\pm i$: $\psi=c_1e^{\alpha y}\cos \beta y+c_2e^{\alpha y}\sin \beta y,\,
\tau=\tilde{c_1}e^{\alpha y}\cos \beta y+\tilde{c_2}e^{\alpha y}\sin \beta y,$ where $\alpha=-\frac{l_1}{l_2},\,\beta=\frac{1}{l_2}$;

Case $l=l_1+il_2=\pm i$: $\psi=c_1y\cos y+c_2y\sin y,\,\tau=\tilde{c_1}y\cos y+\tilde{c_2}y\sin y$.
After a linear transformation we can put in both cases $c_1=\tilde{c_2}=1,\,c_2=\tilde{c_1}=0.$ Thus we have proved the following proposition.

\prop{2.6} Any totally non-degenerate CR-manifold with $\g$ of
type V is locally CR-equivalent to one of the following surfaces:
\begin{gather*}
 \mbox{Type Va:}\,\,  v_2=e^{x+\alpha y}\cos \beta y,\,v_3=e^{x+\alpha
y}\sin \beta y,\,\beta > 0,\alpha+\beta i\neq
i.\\
 \mbox{Type Vb:}\,\,  v_2=e^x y \cos y,\,v_3=e^{x} y \sin y.\end{gather*}

 \type{VI} Lie algebras of this type are characterized by the property that they have an abelian 3-dimensional ideal. Hence we have the
 following commuting relations:
 $$[X_i,X_j]=0,1\leq i,j\leq 3,\,[X_i,X_4]=\sum\limits_{j=1}^3 c_{ij}X_j,1\leq i\leq 3,\,c_{ij}\in\RR{}.$$
  Applying proposition 2.1, we can rectify the commuting vector fields $\x,\xx,\xxx$ simultaneously, then we have $\x=\dz,\xx=\dww,
  \xxx=\dwww$. From the commuting relations we can conclude now that $\xxxx$ is an affine vector field with a real linear part and
  hence (since $M$ is invariant under $\x,\xx,\xxx$) \it $M$ is a tube over a locally affinely homogeneous curve in $\RR{3}$. \rm
  All possible actions of affine 1-dimensional transformation groups in $\RR{3}$ were classified in \cite{shirokov}.
  It is not difficult to obtain from that classification the affine classification of locally affinely homogeneous curves in $\RR{3}$ and hence   the real-affine classification of the  corresponding tubes in $\CC{3}$.
  Rejecting the totally degenerate surfaces, we get the following proposition.

\prop{2.7}  Any totally non-degenerate CR-manifold with $\g$ of
type VI is locally CR-equivalent to one of the following pairwise
locally affinely non-equivalent affinely homogeneous tube surfaces:
\begin{gather*}
{\mbox{Type}\,\, VIa:}\,\, v_2=y^{\alpha},\,v_3=y^{\beta},\,1<\alpha<\beta.\\
{\mbox{Type}\,\, VIb:}\,\, v_2=e^{ay},\,v_3=e^{y},\,-1\leq a<1.\\
{\mbox{Type}\,\, VIc:}\,\, v_2=y\ln y,\,v_3=y^{\alpha},\,\alpha\neq 0;1.\\
{\mbox{Type}\,\, VId:}\,\, v_2=ye^{y},\,v_3=e^{y}.\\
{\mbox{Type}\,\, VIe:}\,\, v_2=y^{2},\,v_3=e^{y}.\\
{\mbox{Type}\,\, VIf:}\,\, v_2=y\ln^{2}y,\,v_3=y\ln y.\\
{\mbox{Type}\,\, VIg:}\,\, v_2=e^{y}\cos\beta y,\,v_3=e^{y}\sin\beta y,\,\beta>0.\\
{\mbox{Type}\,\, VIh:}\,\,v_2=y^{\alpha}\cos(\beta \ln y),\,v_3=y^{\alpha}\sin(\beta \ln y),\,\beta>0.\\
{\mbox{Type}\,\, VIi:}\,\, v_2=\cos y,\,v_3=\sin y. \end{gather*}

\bf\mbox{Types}\,VII-VIII. \rm Lie algebras of these types have the following commuting
relations:
\begin{gather*}
\notag [\x,\xx]=\x,\,[\x,\xxx]=2\xx,\,[\xx,\xxx]=\xxx,\\
[\x,\xxxx]=[\xx,\xxxx]=[\xxx,\xxxx]=0\,\,\,(\mbox{Type VII}),\\
\notag [\x,\xx]=\xxx,\,[\x,\xxx]=-\xx,\,[\xx,\xxx]=\x,\\
[\x,\xxxx]=[\xx,\xxxx]=[\xxx,\xxxx]=0\,\,\,(\mbox{Type VIII}).\end{gather*}

Consider firstly an algebra of type VIII. It contains the subalgebra $\mathfrak{a}=\mbox{span}\{\x,\xx,\xxx\}$, isomorphic to
$\mathfrak{so}_3(\RR{})$. By proposition 2.1, the values of the vector fields $\x,\xx,\xxx$ at $p$ are linearly independent over $\CC{}$, which implies that the orbit of the natural local action of the complexified algebra $\mathfrak{a}^{\CC{}}$ at $p$ is an open set in $\CC{3}$. Hence there exists a local biholomorphic mapping from the Lie group $SO_3(\CC{})$ to a neighborhood of $p$ such that
$\mathfrak{a}^{\CC{}}$ is the image of the tangent algebra of the Lie group $SO_3(\CC{})$ under this mapping. Hence after a
local holomorphic coordinate change we may suppose that $\mathfrak{a}^{\CC{}}$ is the algebra of left-invariant vector fields on
$SO_3(\CC{})$ and that $\xxxx$ is a vector field on $SO_3(\CC{})$, commuting with this algebra. Considering the flows of vector
fields from  $\mathfrak{a}^{\CC{}}$ and of the vector field $\xxxx$, we conclude that these flows commute, which implies that
any transformation from the flow of $\xxxx$ commute with the standard left multiplications in $SO_3(\CC{})$. If $\varphi$ is one of
these transformations and $x,g\in SO_3(\CC{})$, then we get $\varphi(g\cdot x)=g\cdot\varphi(x)$. Putting $x=e$, we get
$\varphi(g)=g\varphi(e)$, which means that $\xxxx$ generates a one-parametric subgroup of right multiplications and hence is a
right-invariant vector field.

Now the orbit of the given algebra at the point $\mbox{Id}$ (corresponding to the original point $p$) is described as follows. Consider $\mathfrak{a}$ as a real subalgebra in the matrix Lie algebra $\mathfrak{so}_3(\CC{})$. Since the values of the vector fields $\x,\xx,\xxx$ at $p$ are linearly independent over $\CC{}$, this subalgebra is a totally-real subspace. Let $A_1,A_2,A_3$ be a basis of this subspace. Then orbit of the action of the real Lie subgroup, corresponding to $\mathfrak{a}$, is given as $e^{x_1A_1}\cdot e^{x_2A_2}\cdot e^{x_3A_3}\cdot Z,\,Z\in SO_3(\CC{}),x_j\in\RR{}$. The orbit of the action of $\xxxx$ is given as $Z\cdot e^{-Bt},\,Z\in SO_3(\CC{}),t\in\RR{}$ for some matrix $B$ from $\mathfrak{so}_3(\CC{})$. Since $\mathfrak{a}$ is totally real, the matrix $B$ can be presented as $B_1+iB_2,\,B_1,B_2\in\mathfrak{a}$. If $\tilde{\xxxx}$ is the right-invariant vector field, corresponding to the matrix $iB_2$, then (since $B_1\in\mathfrak{a}$) at each point the vector fields $\x,\xx,\xxx,\xxxx$ and $\x,\xx,\xxx,\tilde{\xxxx}$ span the same 4-dimensional real linear space. This observation allows us to put $B_1=0$. Now we choose a basis in $\mathfrak{a}$ in such a way that $A_3=-B_2$. Then we finally obtain that the desired orbit looks as follows:
$$e^{x_1A_1}\cdot e^{x_2A_2}\cdot e^{x_3A_3}\cdot e^{itA_3},\,x_j,t\in\RR{}.$$ Consider now the mapping $F: \CC{3}\lr SO_3(\CC{})$, given as $F(z_1,z_2,z_3)=e^{z_1A_1}\cdot e^{z_2A_2}\cdot e^{z_3A_3}$. It is biholomorphic at the origin, since
$F_{z_j}(0)=A_j$ and $A_j$ are linearly independent over $\CC{}$,
and the orbit turns out to be the image of the 4-plane $\{\im
z_1=\im z_2=0\}$ under $F$. Hence  the orbit is degenerate, which
is a contradiction. In the same way we obtain that all $M$ with
$\g$ of type VII are degenerate (in that case the subalgebra,
spanned by $\x,\xx,\xxx$ is isomorphic to
$\mathfrak{so}_{2,1}(\RR{})$). Thus we have proved the following
proposition.

\prop{2.8} Any CR-manifold $M$ with $\g$ of types VII-VIII is
degenerate. \rm

We resume this chapter by formulating the following partial classification theorem.

\theor{2.9} Any totally non-degenerate locally homogeneous
4-dimensional CR-manifold in $\CC{3}$ is locally CR-equivalent to
one of the homogeneous surfaces
Ia,Ib,Ic,II,IIIa,IIIb,IV,Va,Vb,VIa-VIi.\rm

\remark{2.10} Realizations of low-dimensional real Lie algebras as
algebras of vector fields in a real linear space were considered
in many papers. For example, realizations of 4-dimensional real
Lie algebras as algebras of vector fields in $\RR{3}$ were
considered in \cite{nester1} and some formulas, obtained in the
present section, are presented in \cite{nester1}, but the direct
application of the results of \cite{nester1} is impossible in our
case since the situation of a real algebra of holomorphic vector
fields in a complex space gives some restrictions on the possible
realizations (as, for example, proposition 2.1 shows) as well as
some new possibilities, as the above examples show.

\section{The classification}

In this section we specify the partial classification theorem
2.9. This finally allows us to prove the Main Theorem. More
precisely, we study the local CR-equivalence relations among the
surfaces Ia-VIi. To do so, we firstly give the following
definition.

\defin{3.1} A totally non-degenerate locally homogeneous 4-dimensional CR-manifold $M$ is called \it spherical, \rm if
at some point (and hence at each point) it is locally CR-equivalent to the cubic $C$. Otherwise $M$ is called \it non-spherical. \rm

The term "spherical" is used in analogue with the case of a hypersurface in $\CC{2}$, where the 3-dimensional sphere is the model surface for the class of Levi non-degenerate hypersurfaces \cite{poincare}.
Using the trichotomy for 4-dimensional locally homogeneous CR-submanifolds in $\CC{3}$ (see introduction), we get the
following proposition.

\prop{3.2} Two non-spherical 4-dimensional totally non-degenerate locally homogeneous CR-submanifold $M,M'$ are locally CR-equivalent if and only if there exists a
 biholomorphic mapping $F$, translating a point $p\in M$ to a point $p'\in M'$ and (locally) translating
the algebra $\g=\mathfrak{aut}\,M_p$ into the algebra $\mathfrak{g}(M')=\mathfrak{aut}\,M'_{p'}$. In particular, if $M,M'$ are CR-equivalent,
then $\g,\mathfrak{g}(M')$ are isomorphic as Lie algebras.\rm

It follows from the above proposition that for non-spherical manifolds the local CR-equivalence problem can be reduced to the
biholomorphic equivalence problem for vector field algebras, providing the homogeneity of the manifolds. Hence it is important
now to find out what surfaces in the extended list Ia-VIi are spherical.
We firstly note that the sphericity of VIa with $\alpha=2,\beta=3$ follows from \cite{C3}.
The sphericity of Ic for $\alpha=2$ can be verified from the previous fact by applying the binomial formula for $(x+iy)^3$.
To do the sphericity check for the other  surfaces, we refer to the sphericity criterion, formulated in
\cite{bes2}. According to this criterion, a 4-dimensional totally non-degenerate locally homogeneous CR-manifold $M$ is spherical
if and only if in some local coordinates $(z,w_2,w_3)$ it can be presented as $$v_2=|z|^2+O(6),\,\,v_3=2|z|^2\re z+O(7),$$
where the variables are assigned the weights  $[z]=1,[w_2]=2,[w_3]=3$ and $O(6),O(7)$ are terms of weights grater than 6 and 7
correspondingly (this is some analogue for the sphericity criterion for a hypersurface in $\CC{N}$, see \cite{cartan},\cite{chern}). To apply this criterion to the above list of homogeneous surfaces, we consider all possible holomorphic transformations, preserving the origin, and present them as  \begin{gather*}z \mapsto f_1+ \dots + f_n+O(n+1), \quad w_2 \mapsto g_1+ \dots + g_{n+1}+O(n+2),\\
w_3 \mapsto h_1+ \dots + h_{n+2}+O(n+3),\end{gather*} where $f_j,g_j,h_j$ are polynomials of weight $j$. We call the collection of all $f_j,g_{k},h_{l}$
for $j\leq n,\,k\leq n+1,\,l\leq n+2$ \it the $(n,n+1,n+2)$-jet \rm of the transformation. We also present the manifold as
$$v_2=|z|^2+\sum\limits_{j=3}^{\infty}F_j,\,v_3=2|z|^2\re z+\sum\limits_{j=4}^{\infty}G_j,$$ where $F_j,G_j$ are polynomials of
weight $j$ and call the collection of all $F_j,G_{j+1}$ for $3\leq j\leq m$ \it the $(m,m+1)$ jet \rm of $M$. Then we note the following. Given a mapping of a manifold
\begin{gather*}v_2=|z|^2+F_3+ \dots + F_m+O(m+1),\\ v_3= 2|z|^2\re z+G_4+\dots+G_{m+1}+O(m+2)\end{gather*}
to a manifold of the same kind \begin{gather*}v_2=|z|^2+\hat{F}_3+ \dots + \hat{F}_m+O(m+1),\\ v_3=2 |z|^2\re z+\hat{G}_4+\dots+\hat{G}_{m+1}+O(m+2),\end{gather*}
preserving the origin, for a fixed $(m,m+1)$-jet of the first manifold the $(m,m+1)$-jet of the second manifold depends only on the  $(m-1,m,m+1)$-jet of the mapping.
Then it is clear that the vanishing condition for the  $(5,6)$-jet of the manifold (which is equivalent to the sphericity) is a condition on the  $(4,5,6)$-jet of the mapping. This condition is
a system of equations on the coefficients of the  $(4,5,6)$-jet of the mapping.

Now we describe the process of the sphericity inspection for a totally non-degenerate surface.

{\bf Step 0}. We present the equations of the surface as
$$v_2=|z|^2+F_3+ F_4 + F_5+O(6), \quad  v_3= 2|z|^2\re z+G_4+G_5+G_6+O(7).$$

{\bf Step 1}.  We write the condition on the coefficients of a mapping
$$z \mapsto z + f_2+O(3), \quad w_2 \mapsto w_2 + g_3+O(4), \quad
w_3 \mapsto w_3+ h_4+O(5),$$
which maps the original surface onto a surface
$$v_2=|z|^2+O(4), \quad  v_3= 2|z|^2\re z+O(5).$$

This condition is a system of 17 real equations on  18 real variables.  This system always has a solution $(f_2,g_3,h_4)$.

{\bf Step 2}.  We write the condition on the coefficients of a mapping
$$z \mapsto z + f_2+f_3+O(4), \, w_2 \mapsto w_2 + g_3+g_4+O(5),\,
w_3 \mapsto w_3+ h_4+h_5+O(6),$$
which maps the surface, obtained in step 1, onto a surface
$$v_2=|z|^2+O(5), \quad  v_3= 2|z|^2\re z+O(6).$$
This condition is a system of 26 real equations on 24 real variables,  which might not have any solution. If this system has no solution, then we conclude that
the surface is not spherical. Otherwise we get a solution $(f_3,g_4,h_5)$ and go to step 3.

{\bf Step 3}.  We write the condition on the coefficients of a mapping
\begin{gather*}z \mapsto z + f_2+f_3+f_4+O(5), \quad w_2 \mapsto w_2 + g_3+g_4+g_5+O(6), \\
w_3 \mapsto w_3+ h_4+h_5+h_6+O(7),\end{gather*}
which maps the surface, obtained in step 2, onto a surface
$$v_2=|z|^2+O(6), \quad  v_3= 2|z|^2\re z+O(7).$$
This condition is a system of 39 real equations on 32 real
variables,  which may have no solution. If this system has no
solution, then we conclude that the surface is not spherical,
otherwise it is spherical.

To solve the systems of equations for the homogeneous surfaces Ia-Vb we used Maple 6 package (see \cite{appendix} for the details of the computations).
For the tube case VI it is possible to apply simpler arguments. We resume the results of our computations in the following proposition.

\prop{3.3} The following homogeneous surfaces from the list Ia-VIi
are spherical: Ic for $\alpha=2$ and VIa for
$\alpha=2,\beta=3$. All other surfaces from the list Ia-VIi are
non-spherical.

\doc To prove the proposition for the tube surfaces VIa-VIi we note the following: the subalgebra
$\mbox{span}\{\xxx,\xx,\xp\}$ is the unique abelian 3-dimensional subalgebra in the 5-dimensional infinitesimal automorphism
algebra of the cubic $C$ (it can be easily verified from the commuting relations in the algebra, see \cite{C3} for the details and the notations).
As it was shown in \cite{C3}, the cubic $C$ is polinomially equivalent to the tube surface $VIa$ for $\alpha=2,\beta=3$
(we denote this surface by $\tilde{C}$) and we conclude that $\mbox{span}\{\dz,\dww,\dwww\}$ is the unique 3-dimensional
abelian subalgebra
in $\mathfrak{aut}\tilde{C}$. Hence if a tube surface $M$ from the list VIa-VIi is locally biholomorphically equivalent
to the cubic $C$, we get a biholomorphic mapping $F$, which maps $\g$ to a 4-dimensional subalgebra of the 5-dimensional algebra
$\mathfrak{aut}\tilde{C}$. In particular, the abelian 3-algebra, spanned by $\x,\xx,\xxx$, is mapped to an abelian 3-dimensional
subalgebra of $\mathfrak{aut}\tilde{C}$. Since such subalgebra is unique, we conclude that $F$ maps the three coordinate vector field
$\dz,\dww,\dwww$ to their linear combinations, which implies that $F$ is in fact a linear mapping. Hence all the tube surfaces VIa-VIi
except $\tilde{C}$ are locally biholomorphically non-equivalent to $C$, as required. \qed

It remains now to prove that the non-spherical homogeneous
surfaces Ia-VIi are pairwise locally CR-inequivalent. Since
homogeneous surfaces of different types Ia-VIi correspond to
non-isomorphic Lie algebras (except the cases Va,Vb,
corresponding to the same algebra V), we just need to prove, using
proposition 3.2, that two vector field algebras, providing the
homogeneity of two non-spherical surfaces of the same type Ia-VIi,
can not be mapped onto each other by a local biholomorphic mapping
and also to prove that two vector field algebras, providing the
homogeneity of a type Va surface and a type Vb surface correspondingly, can not be
mapped onto each other by a local biholomorphic mapping. In what
follows
$$\Phi:\,\,\CC{3}_{z,w_2,w_3}\mapsto\CC{3}_{\xi,\eta_2,\eta_3}$$
denotes a local biholomorphic mapping, which maps a germ of a
homogeneous surface $M$ at $p$ onto a germ of a homogeneous
surface $M'$ at $p'$. $\x,\xx,\xxx,\xxxx$ and $Y_1,Y_2,Y_3,Y_4$
denote the basis of the vector field algebras $\g$ and
$\mathfrak{g}(M')$ correspondingly. Now we consider different
cases.

$\bf Ia\mapsto Ia$. \rm In this case we may suppose
$p=p'=(0,0,i)$.  From section II we have
$\x=\dww,\,\xx=\dwww,\,\xxx=\dz+w_3\dww,\,\xxxx=(i-\gamma)\dz+w_2\dww+w_3\dwww$
and the same for $Y_j$ (with a parameter $\gamma'$). The vector
field algebras are of type I with $q=0$. $\mbox{Span}\{\x,\xx\}$
is the commutant and hence is invariant under $\Phi$. $\xxxx$ is
the unique element, modulo the commutant, for which the
corresponding adjoint operator is identical on the commutant.
$\xxx$ is the unique element, up to a scalar and modulo the
commutant, for which the corresponding adjoint operator has zero
eigenvalues on the commutant. $\mbox{Span}\{\x\}$ is the kernel of
$\mbox{ad}_{\xxx}$. From all the above invariant descriptions we
conclude that $\Phi$ maps $\x$ to $aY_1$, $\xx$ to $bY_1+eY_2$,
$\xxx$ to $mY_1+nY_2+kY_3$, $\xxxx$ to $pY_1+sY_2+Y_4$, where
$a,b,e,m,n,k,p,s\in\RR{},\,a,e,k\neq 0$. The first two conditions
imply $\xi=F(z),\,\eta_2=G(z)+aw_2+bw_3,\,\eta_3=H(z)+ew_3$. The
third one implies $F_z=k$, the fourth one implies
$(i-\gamma)F_z=i-\gamma'$. Hence we get $k=1,\gamma=\gamma'$ and
conclude that $M=M'$.

$\bf Ib\mapsto Ib$. \rm In this case we may suppose
$p=p'=(i,0,i)$.  From section II we have
$\x=\dww,\,\xx=\dwww,\,\xxx=\dz+w_3\dww,\,\xxxx=-z\dz-i\gamma\dww+w_3\dwww$
and the same for $Y_j$ (with a parameter $\gamma'$). The vector
field algebras are of type I with $q=-1$.
$\mbox{Span}\{\x,\xx,\xxx\}$ is the commutant and hence is
invariant under $\Phi$. $\xxxx$ is the unique element, up to a
sign and modulo the commutant, for which the corresponding adjoint
operator has eigenvalues $\{0;1;-1\}$ on the commutant.
$\mbox{Span}\{\x\}$ is the second commutant. From the above
invariant descriptions we conclude that $\Phi$ maps $\x$ to
$aY_1$,  $\xx$ to  $bY_1+eY_2+kY_3,$ $X_3$ to $lY_1+mY_2+nY_3$,
$\xxxx$ to $pY_1+sY_2+rY_3 \pm Y_4$. In case we have plus, the
commuting relations imply $k=m=0$. Then from the first two
conditions we get
$\xi=F(z),\,\eta_2=G(z)+aw_2+bw_3,\,\eta_3=H(z)+ew_3$. Comparing
the $\frac{\partial}{\partial \xi},\,\frac{\partial}{\partial
\eta_2}$ and $\frac{\partial}{\partial \eta_3}$ - coefficients for
the third condition, we get $F_z=n,\,G_z+aw_3=l+nH+new_3,\,H_z=m$
and hence $F=nz-ni+i,\,a=ne,\,G_z=nH+l,H_z=0$ (because $F(i)=i$).
Comparing the $\frac{\partial}{\partial
\xi},\,\frac{\partial}{\partial \eta_2}$  - coefficients for the
fourth condition, we get $-nz=-nz+ni-i+r,\,-zG_z-i\gamma
a+bw_3=p+rH+rew_3-i\gamma',$ which implies
$n=1,r=0,\,G_z=0,\,\gamma'=\gamma a$. This finally gives us
$H=-l\in\RR{}$ and hence $e=1$ (because $\Phi(p)=p'$). Now we
conclude that $a=ne=1$, $\gamma'=\gamma a=\gamma$ and $M=M'$, as
required. In case we have minus in the fourth condition, we note
that $M$ has the following polynomial automorphism $\sigma$,
preserving the point $p$:
$$z\mapsto w_3,\,w_2\mapsto -w_2+zw_3+1,\,w_3\mapsto z.$$
Under this transformation \begin{gather}
\x\mapsto -\x,\,\xx\mapsto\xxx,\,\xxx\mapsto\xx,\,\xxxx\mapsto -\xxxx.\end{gather} Composing $\Phi$ with $\sigma$,
we get a mapping $\tilde\Phi$ from $M$ onto $M'$ with a plus in the fourth condition and hence conclude that $\gamma'=\gamma$ and $M'=M$.

$\bf Ic\mapsto Ic$. \rm In this case we may suppose that the
alpha's are the same for $M$ and $M'$ (since different alpha's
correspond to different $q$ and hence to non-isomorphic Lie
algebras). Also we suppose that
$p=(i,i\gamma,i),\,p'=(i,i\gamma',i)$. The vector field algebras
are of type I with $q\in(-1,1),q\neq 0$. From section II we have
$\x=\dww,\,\xx=\dwww,\,\xxx=\dz+w_3\dww,\,\xxxx=qz\dz+(q+1)w_2\dww+w_3\dwww$
and the same for $Y_j$.  $\mbox{Span}\{\x,\xx,\xxx\}$ is the
commutant.  $\mbox{Span}\{\x\}$ is the second commutant. $\xxxx$ is
the unique element, modulo the commutant, for which the
corresponding adjoint operator has the eigenvalue $q+1$ on the
second commutant.From the above invariant descriptions we conclude
that $\Phi$ maps $\x$ to $aY_1$,  $\xx$ to  $bY_1+eY_2+kY_3,$
$X_3$ to $lY_1+mY_2+nY_3$, $\xxxx$ to $pY_1+sY_2+rY_3+Y_4$. The
commuting relations imply $k=m=0$. Then from the first two
conditions we get
$\xi=F(z),\,\eta_2=G(z)+aw_2+bw_3,\,\eta_3=H(z)+ew_3$. Comparing
the $\frac{\partial}{\partial \xi},\,\frac{\partial}{\partial
\eta_2}$ and $\frac{\partial}{\partial \eta_3}$ - coefficients for
the third condition, we get $F_z=n,\,G_z+aw_3=l+nH+new_3,\,H_z=m$
and hence $a=ne,\,G_z=nH+l,H_z=0$.  In particular, $F=nz-ni+i$
(from $F(i)=i$) and $G$ is linear. Comparing the
$\frac{\partial}{\partial \xi}$ - coefficients for the fourth
condition, we get $qnz=q(nz-ni+i)+r$,which implies $n=1,r=0$.
Comparing the $\frac{\partial}{\partial \eta_2}$ - coefficients
for the fourth condition, we get
$zG_z+a(q+1)w_2+bw_3=p+rH+rew_3+a(q+1)w_2+(q+1)w_3+(q+1)G+p$,
which implies $b=b(q+1)$ and hence $b=0$, and also
$zG_z=\frac{q+1}{q}G$ and $G=Az^{\frac{q+1}{q}}$. Since $G$ is
linear and $q\neq -1$, we get $G=0$. Then $H=-l\in\RR{}$ and hence
$e=1$ (because $\Phi(p)=p'$). Now we conclude that $a=ne=1$, and
hence $\eta_2=w_2$, which implies $i\gamma=i\gamma'$ and $M=M'$,
as required.

$\bf II\mapsto II$. \rm In this case
$p=(i,i\gamma,0),\,p'=(i,i\gamma',0)$. The vector
field algebras are of type II. From section II we have
$\x=\dww,\,\xx=\dwww,\,\xxx=\dz+w_3\dww,\,\xxxx=z\dz+(2w_2+\frac{1}{2}z^2)\dww+(z+w_3)\dwww$
and the same for $Y_j$.  $\mbox{Span}\{\x,\xx,\xxx\}$
is the commutant.  $\mbox{Span}\{\x\}$ is the
second commutant. $\xxxx$ is
unique element, modulo the commutant, for which the
corresponding adjoint operator has the eigenvalue $2$ on the second commutant.From the above invariant descriptions we
conclude that $\Phi$ maps $\x$ to $aY_1$,  $\xx$ to  $bY_1+eY_2+kY_3,$
$X_3$ to $lY_1+mY_2+nY_3$, $\xxxx$ to $pY_1+sY_2+rY_3+Y_4$.
The commuting relations imply $k=0,n=e,a=e^2$. Then in the same way as in the pervious case
from the first two conditions we get $\xi=F(z),\,\eta_2=G(z)+e^2w_2+bw_3,\,\eta_3=H(z)+ew_3$, and from the third one
$F_z=nz-ni+i,\,G_z=eH+l,H_z=m$.  Comparing the $\frac{\partial}{\partial \xi}$ - coefficients for the fourth condition,
we get $nz=nz-ni+i+r$,which implies $n=1,r=0$.  Comparing the $\frac{\partial}{\partial \eta_2}$ - coefficients for the fourth condition, we get $zG_z+(2w_2+\frac{1}{2}z^2)e^2+b(z+w_3)=p+2e^2w_2+2bw_3+2G$, which implies $e^2=1,b=0,p=0,\,zG_z=2G$ and hence  $G=Az^2$. Finally, comparing the $\frac{\partial}{\partial \eta_3}$ - coefficients, we get $zH_z+(z+w_3)e=z+H+ew_3+q$, which
implies $zH_z=H+q,\,H=mz-q$. Now from $\Phi(p)=p'$ we get $mi-q=0$ and hence $m=q=0,\,H=0,\,G=0$ (because $G=Az^2$ and
$G_z=l+eH$). Applying $\Phi(p)=p'$ again, we get $i\gamma=i\gamma'$ and $M=M'$, as required.

$\bf IIIa\mapsto IIIa$. \rm In this case we may suppose
$p=p'=(0,0,i)$. The vector field algebras are of type III with
$q=0$. $\mbox{Span}\{\x,\xx,\xxx\}$ is the commutant,
$\mbox{span}\{\x\}$ is the second commutant, $\xxxx$ is the unique
element, up to a sign and modulo the commutant, for which the
corresponding adjoint operator has eigenvalues $\{0,i,-i\}$ on the
commutant. From the above invariant descriptions we conclude that
$\Phi$ maps $\x$ to $aY_1$,  $\xx$ to  $bY_1+cY_2+dY_3,$ $X_3$ to
$eY_1+kY_2+lY_3$, $\xxxx$ to $mY_1+nY_2+pY_3 \pm Y_4$. From
section II we have $\x=\dww,\,\xx=\dwww,\,\xxx=\dz+w_3\dww,\,
\xxxx=-w_3\dz+(\frac{1}{2}z^2-\frac{1}{2}w_3^2-i\gamma)\dww+z\dwww$
and the same for $Y_j$ (with a parameter $\gamma'$). It is
convenient now to make the variable change (2) and thus to get the
affine vector fields
$\x=\dww,\,\xx=\dwww-w_2\dz,\,\xxx=\dz+w_3\dww,\,\xxxx=-w_3\dz-2i\gamma\dww+z\dwww$,
and the same for $Y_j$ (with a parameter $\gamma'$). The commuting
relations now imply $l=c,\,k=-d,\,c^2+d^2=a$. Then from the first
condition we get $\xi_{w_2}=(\eta_3)_{w_2}=0,\,(\eta_2)_{w_2}=a$.
The third and the second condition imply
$\xi_{w_3}=d,\,\xi_z=c,\,(\eta_3)_z=-d,\,(\eta_3)_{w_3}=c$ and
hence $$\xi=cz+dw_3+s,\,\eta_3=-dz+cw_3+t.$$ Also we get
$(\eta_2)_z=e+sd+ct,\,(\eta_2)_{w_3}=b-cs+dt$. From $\Phi(p)=p'$
we get $s=-di,\,t=i-ci$.   In case we have plus in the fourth
condition, we compare the $\frac{\partial}{\partial \xi}$ and
$\frac{\partial}{\partial \eta_3}$ - coefficients for the fourth
condition and get $p=n=s=t=0$. Hence $c=1,\,d=0,\,a=1$. Thus we
have $\eta_2=w_2+ez+bw_3+h$. Comparing now the
$\frac{\partial}{\partial \eta_2}$ - coefficients for the fourth
condition, we get $b=e=0,\,\gamma=\gamma'$, so $M=M'$. In case we
have minus in the fourth condition, we firstly apply the following
authomorphism $\varepsilon$, preserving $M$ and the fixed point
$p$:
$$z\mapsto -z,\,w_2\mapsto -w_2,\,w_3\mapsto w_3.$$ As a result we get
$\x\mapsto -\x,\,\xx\mapsto \xx,\,\xxx\mapsto -\xxx,\,\xxxx\mapsto -\xxxx.$ Composing $\Phi$ with $\varepsilon$, we get a
biholomorphic mapping $\Phi'$ of a germ of $M$ at $p$ onto a germ of $M'$ at $p'$ with a plus in the fourth condition and hence
get $M=M'$, as required.

$\bf IIIb\mapsto IIIb$. \rm In this case we may suppose that the
$q$-parameters are the same for $M$ and $M'$ (since different $q$
correspond to non-isomorphic Lie algebras). Also we suppose that
$p=(i,i\gamma,0),\,p'=(i,i\gamma',0)$. The vector field algebras
are of type III with $q>0$. $\mbox{Span}\{\x,\xx,\xxx\}$ is the
commutant,  $\mbox{span}\{\x\}$ is the second commutant, $\xxxx$
is the unique element, modulo the commutant, for which the
corresponding adjoint operator has the eigenvalue $2q$ on the
second commutant. From the above invariant descriptions we
conclude that $\Phi$ maps $\x$ to $aY_1$,  $\xx$ to
$bY_1+cY_2+dY_3,$ $X_3$ to $eY_1+kY_2+lY_3$, $\xxxx$ to
$mY_1+nY_2+pY_3+Y_4$. From section II we have
$\x=\dww,\,\xx=\dwww,\,\xxx=\dz+w_3\dww,\,\xxxx=(qz-w_3)\dz+(2qw_2+\frac{1}{2}z^2-\frac{1}{2}w_3^2)\dww+(z+qw_3)\dwww$
and the same for $Y_j$. It is convenient now to make the variable
change (2) and thus to get the affine vector fields
$\x=\dww,\,\xx=\dwww-w_2\dz,\,\xxx=\dz+w_3\dww,\,\xxxx=(qz-w_3)\dz+2qw_2\dww+(z+qw_3)\dwww$
(and the same for $Y_j$). The commuting relations now imply
$l=c,\,k=-d,\,c^2+d^2=a$. Then from the first condition we get
$\xi_{w_2}=(\eta_3)_{w_2}=0,\,(\eta_2)_{w_2}=a$. The third and the
second condition imply
$\xi_{w_3}=d,\,\xi_z=c,\,(\eta_3)_z=-d,\,(\eta_3)_{w_3}=c$ and
hence $$\xi=cz+dw_3+s,\,\eta_3=-dz+cw_3+t.$$ Also we get
$(\eta_2)_z=e+sd+ct,\,(\eta_2)_{w_3}=b-cs+dt$. From $\Phi(p)=p'$
we get $s=i-ci,\,t=di$.   Comparing the $\frac{\partial}{\partial
\xi}$ and $\frac{\partial}{\partial \eta_3}$ - coefficients for
the fourth condition, it is not difficult to obtain $p=n=s=t=0$
and hence $c=1,\,d=0,\,a=1$. Thus we have $\eta_2=w_2+ez+bw_3+h$.
Comparing now the  $\frac{\partial}{\partial \eta_2}$ -
coefficients for the fourth condition, we get $b=e=0,\,m+2qh=0$,
which implies $\im h=0$ and from $\Phi(p)=p'$ we now get
$i\gamma=i\gamma'$ and $M=M'$, as required.

$\bf IV\mapsto IV$. \rm In this case we may suppose
$p=p'=(0,i,i)$.  From section II we have
$\x=\dww,\,\xx=\dwww,\,\xxx=\dz+w_3\dwww,\,\xxxx=(i-\delta)\dz+w_2\dww$
and the same for $Y_j$ (with a parameter $\delta'$). The vector
field algebras are of type IV. $\mbox{Span}\{\x,\xx\}$ is the
commutant. $\xxx$ and $\xxxx$ are the unique elements, up to a
permutation and modulo the commutant, for which the corresponding
adjoint operators have a collection of eigenvalues $\{0;1\}$ on
the commutant. From the above invariant descriptions we conclude
that $\Phi$ maps $\x$ to $aY_1+pY_2$, $\xx$ to $bY_1+qY_2$, and
also (first case) $\xxx$ to $mY_1+nY_2+Y_3$, $\xxxx$ to
$kY_1+lY_2+Y_4$, or (second case) $\xxx$ to $mY_1+nY_2+Y_4$,
$\xxxx$ to $kY_1+lY_2+Y_3$. The first two conditions in both cases
imply $\xi=F(z),\,\eta_2=G(z)+aw_2+bw_3,\,\eta_3=H(z)+pw_2+qw_3$.
Comparing now the $\frac{\partial}{\partial \xi}$ - coefficients
for the third and the fourth conditions in the first case, we get
$F_z=1$ and $(i-\delta)F_z=i-\delta'$, which implies
$\delta=\delta'$ and $M=M'$, as required. Comparing the
$\frac{\partial}{\partial \xi}$ - coefficients for the third and
the fourth conditions in the second case, we get $F_z=i-\delta$
and $(i-\delta)F_z=1$, which implies $(i-\delta)^2=1$, which is
impossible since $\delta\in\RR{}$, so the second case can not
occur and finally $M=M'$.

$\bf V\mapsto V$. \rm We consider a mapping between germs of
arbitrary surfaces with $\g$ of type V. In this case from section
II we have $\x=\dww,\,\xx=\dwww,\,\xxx=\dz+w_2\dww+w_3\dwww$. For
type Va we have $X_4=\frac{i-\alpha}{\beta}\dz-w_3\dww+w_2\dwww,$
for type Vb we have $X_4=i\dz+(ie^z-w_3)\dww+(w_2+e^z)\dwww.$
$\mbox{Span}\{\x,\xx\}$ is the commutant. $\xxx$ is the unique
element, modulo the commutant, for which the corresponding adjoint
operator is identical on the commutant. $\xxxx$ is the unique
element, up to a sign and modulo the commutant, for which the
corresponding adjoint operator has eigenvalues $\{i;-i\}$. From
the above invariant descriptions we conclude that $\Phi$ maps $\x$
to $aY_1+pY_2$, $\xx$ to $bY_1+qY_2$, $\xxx$ to $mY_1+nY_2+Y_3$,
and also (first case) $\xxxx$ to $kY_1+lY_2+Y_4$, or (second case)
$\xxxx$ to $kY_1+lY_2-Y_4$. The first two conditions in both cases
imply $\xi=F(z),\,\eta_2=G(z)+aw_2+bw_3,\,\eta_3=H(z)+pw_2+qw_3$.
Comparing now the $\frac{\partial}{\partial \xi}$ - coefficients
for the third and the fourth conditions in the first case, we get
$F_z=1$, and also:
$\frac{i-\alpha}{\beta}F_z=\frac{i-\alpha'}{\beta'}$ ($Va\mapsto
Va$), or  $\frac{i-\alpha}{\beta}F_z=i$ ($Va\mapsto Vb$), which
implies $\alpha=\alpha',\beta=\beta'$ and $M=M'$ for $Va\mapsto
Va$, or $\alpha=0,\beta=1$ for $Va\mapsto Vb$, which is a
contradiction. Comparing the $\frac{\partial}{\partial \xi}$ -
coefficients for the third and the fourth conditions in the second
case, we get $F_z=1$, and also:
$\frac{i-\alpha}{\beta}F_z=-\frac{i-\alpha'}{\beta'}$ ($Va\mapsto
Va$), or $\frac{i-\alpha}{\beta}F_z=-i$ ($Va\mapsto Vb$), which
implies $\beta=-\beta'$ for $Va\mapsto Va$, which is a
contradiction since $\beta,\beta'>0$, or $\alpha=0,\beta=-1$ for
$Va\mapsto Vb$, which is also a  contradiction since $\beta>0$, as
required.

$\bf VI\mapsto VI$. \rm We claim that all non-sperical tube
surfaces VIa-VIi are pairwise non-equivalent. To see that, we note
that for all tube surfaces VIa-VIi the real matrix, defining the
linear part of the affine vector field $X_4$, has rank at least 2,
which implies that the centralizer of $X_4$ in $\g$ is not more
than 1-dimensional. Then we can find no 3-dimensional abelian
subalgebra in $\g$ other than the one spanned by $\x,\xx,\xxx$.
Hence if $F$ is a biholomorphic mapping between germs of two
non-spherical tube surfaces $M,M'$ from the list VIa-VIi, then by
proposition 3.2 $F$ maps $\g$ to $\mathfrak{g}(M')$ and the unique
3-dimensional abelian subalgebras are also mapped onto each other.
In the same way as in proposition 3.3 we conclude now that $F$ is
a linear mapping. Hence different tube surfaces from the list
VIa-VIi are pairwise locally CR non-equivalent, as required.

\smallskip

Thus, according to the claims of the main trichotomy, theorem 2.9
and proposition 3.3, the Main Theorem is completely proved. \qed

\smallskip

\noindent{\bf Proof of theorem 1.2}. \rm For all surfaces
(2.1)-(3.19) except (3.4) and (3.9)  the desired claim follows
from the fact that in appropriate coordinates $\g$ consists
actually of affine vector fields (see the above description of
$\g$ for different types). For the exceptional surface (3.4) one
should make the variable change (2). After that we have (see the
case II$\mapsto$II above):
$\x=2\dww,\,\xx=\dwww-z\dww,\,\xxx=\dz+w_3\dww,\,\xxxx=z\dz+2w_2\dww+(z+w_3)\dwww$
and thus the homogeneity of (3.4) is provided by an algebra of
affine vector fields. For the exceptional case (3.9) one should
make the variable change
$z^*=e^z$. Then $\g$ looks as follows (see the case V$\mapsto$V above): \begin{gather*}\x=\dww,\,\xx=\dwww,\,\xxx=z\dz+w_2\dww+w_3\dwww,\\
\xxxx=iz\dz+(iz-w_3)\dww+(z+w_2)\dwww\end{gather*} and hence consists of affine vector fields, as required. \qed

\noindent{\bf Proof of theorem 1.2'}. \rm For each algebra of type
I-V except the type I algebra with $q=1$ we can find some
homogeneous surface in the list (3.1) - (3.9), for which this
algebra is a transitively acting   algebra (as it follows from
section 2), and then to apply theorem 1.2. For the type I algebra
with $q=1$ the necessary realization is given in the case I of
section 2. For type VI algebras one should just put
$\x=\frac{\partial}{\partial z_1},\,\xx=\dzz,\,\xxx=\dzzz$ and
$\xxxx$ to be a linear vector field with the defining matrix
$C=(c_{ij}).$ For a type VIII algebra one should put (say, in a
neighborhood of the point $(0,1,i)\in\CC{3}$) $X_j$ with $j=1,2,3$
to be a linear vector field with a defining matrix from the
standard basis of the Lie algebra $\mathfrak{so}_3(\CC{})$, and
also put $\xxxx$ to be the Euler vector field. For a type VII
algebra the proof is similar. \qed

\noindent{\bf Proof of theorem 1.3}. \rm Following carefully the
above arguments, one can see that in cases Ia-Vb  any mapping of a
germ of a homogeneous surface onto itself is actually identical
except the mapping $\sigma$ for the surfaces of type Ib and
$\varepsilon$ for the surfaces of type IIIa. For tube surfaces
VIa-VIi each mapping of a germ onto itself must be linear (see the
arguments above). It is not difficult to see that only the cubic,
the surfaces VIb for $a=-1$ and the surface VIi have the desired
linear automorphism, preserving a fixed point (the affinely
homogeneous curves, corresponding to tube surfaces VIa, can not be
extended to the origin as affinely homogeneous curves except the
case of the cubic). \qed

\noindent{\bf Proof of theorem 1.4}. \rm As it was mentioned
above, any symmetric CR-manifold is CR-homogeneous. It follows
from theorem 1.3 that among the totally non-degenerate homogeneous
surfaces (2.1) - (3.19) the following ones are symmetric: (2.1),
(3.2), (3.5), (3.12) and (3.19), the desired involutions $s_p$ are
as they are presented in theorem 1.3; for the cubic one should put
$\lambda=-1$. The isometry property for $s_p$ is obvious in cases
(2.1), (3.5), (3.12), (3.19). In the case of (3.2) the isometry
property follows from (6). The symmetry property check in all
cases is straightforward. For the  degenerate surfaces (1.1) -
(1.3) we note that a CR-manifold of kind $M^3\times \RR{1}$, where
$M^3\subset\CC{2}_{z,w_2}$ is a locally homogeneous surface in
$\CC{2}$, $\RR{1}\subset\CC{1}_{w_3}$ is a real line, is
CR-symmetric if and only if $M^3$ is symmetric. It follows from
the E.Cartan's classification theorem that among the hypersurfaces
from his list only the surfaces, corresponding to (2.b) - (2.f),
have non-trivial stability subgroups, which are of
$\mathbb{Z}_2$-type. The symmetry property check for these
surfaces, as well as for a hyperplane and a hypersphere in
$\CC{2}$, is straightforward. This completely proves the theorem.
\qed

\remark{3.4} It is an amazing consequence of theorem 1.4 that any
symmetric totally non-degenerate 4-manifold is associated to a
second order plane curve in the sense that  one of the defining
equations of the manifold can be chosen as an equation of a second
order plane curve.

\remark{3.5} Another possible approach to the description of
homogeneous surfaces is to present a list of all possible normal
forms (see, for example, \cite{ezheas}, \cite{loboda3}). For some
cases this approach is actually realized in \cite{appendix}. It
would be interesting, taking corollary 1.2 into account, to
reformulate the Main Theorem in terms of some affine normal forms,
as it was made, for example, in \cite{ezheas} for affinely
homogeneous hypersurfaces in $\RR{3}$. It would be also
interesting to find the specify of the normal form for the
symmetric surfaces.

\bigskip
\small{\obeylines
 Valery K.Beloshapka
 Department of Mathematics
 The Moscow  State University
 Leninskie Gori, MGU, Moscow, RUSSIA
 E-mail: vkb@strogino.ru
 --------------------------------------------------------------------
 Ilya G.Kossovskiy
 Department of Mathematics
The Australian National University
 Canberra, ACT 0200 AUSTRALIA
 E-mail: ilya.kossovskiy@anu.edu.au }


\begin{thebibliography}{99}

\bibitem{nachin} A.Altomani, C.Medori, M.Nacinovich, "On homogeneous and symmetric CR manifolds", 2009, arXiv:0910.4531.

\bibitem{bouendy}  M.S.Baouendi, P.Ebenfelt, L.P.Rothschild,
"Real Submanifolds in Complex Space and Their Mappings",
Princeton University  Press, Princeton Math. Ser. 47, Princeton,
NJ, 1999.

\bibitem{univers} Beloshapka, V. K. "A universal model for a real submanifold", (Russian)  Mat. Zametki  75  (2004),  no. 4, 507--522;  translation in  Math. Notes  75  (2004),  no. 3-4, 475--488.

\bibitem{small} Beloshapka, V. K. "CR-varieties of the type $(1,2)$ as varieties of ``super-high'' codimension",  Russian J. Math. Phys.  5  (1997),  no. 3, 399--404 (1998).

\bibitem{bes}  Beloshapka, V. K.; Ezhov, V. V.; Shmalts, G. "Holomorphic classification of four-dimensional surfaces in $\Bbb C\sp 3$" (Russian),  Izv. Ross. Akad. Nauk Ser. Mat.  72  (2008),  no. 3, 3--18.

\bibitem{bes2} Beloshapka, V. K.; Ezhov, V. V.; Shmalts, "Canonical Cartan connection and holomorphic invariants on Engel CR manifolds",  Russ. J. Math. Phys.  14  (2007),  no. 2, 121--133.

\bibitem{appendix} Beloshapka, V.K., Kossovskiy, I.G., Appendix to the paper "Classification of homogeneous CR-manifolds in
dimension 4", www.arxiv.org, to be appeared soon.

\bibitem{C3} Beloshapka, V.K., Kossovskiy, I.G., "Homogeneous hypersurfaces in $\CC{3}$, associated with a model
 CR-cubic", Jounal of Geometric Analysis, 2010 (in print); arXiv:0910.0658.

\bibitem{cartan} Cartan, E. "Sur la geometrie pseudo-conforme des hypersurfaces de l'espace de deux variables complexes II". (French)  Ann. Scuola Norm. Sup. Pisa Cl. Sci. (2)  1  (1932),  no. 4, 333--354.

\bibitem{chern} Chern, S. S.; Moser, J. K. "Real hypersurfaces in complex manifolds",  Acta Math.  133  (1974), 219--271.

\bibitem{ezheas} Eastwood, M.; Ezhov, V., "On affine normal forms and a classification of homogeneous surfaces in affine three-space",  Geom. Dedicata  77  (1999),  no. 1, 11--69.

\bibitem{kaup} Fels, G.; Kaup, W., "Classification of Levi degenerate homogeneous CR-manifolds in dimension 5",  Acta Math.  201  (2008),  no. 1, 1--82.

\bibitem{kaup1} Fels, G.; Kaup, W., "CR-manifolds of dimension 5: a Lie algebra approach",  J. Reine Angew. Math.  604  (2007), 47--71.

\bibitem{kaupzaitsev} Kaup, W.; Zaitsev, D., "On symmetric Cauchy-Riemann manifolds",  Adv. Math.  149  (2000),  no. 2, 145--181.


\bibitem{loboda1} Loboda, A. V. "Homogeneous real hypersurfaces in $\CC{3}$ with two-dimensional isotropy groups", (Russian)  Tr. Mat. Inst. Steklova  235  (2001),  Anal. i Geom. Vopr. Kompleks. Analiza, 114--142;  translation in  Proc. Steklov Inst. Math.  2001,  no. 4 (235), 107--135.

\bibitem{loboda2} Loboda, A. V. "On the dimension of a group that acts transitively on hypersurfaces in $C\sp 3$". (Russian)  Funktsional. Anal. i Prilozhen.  33  (1999),  no. 1, 68--71;  translation in  Funct. Anal. Appl.  33  (1999),  no. 1, 58--60

\bibitem{loboda3} Loboda, A. V. "On the determination of a homogeneous strictly pseudoconvex hypersurface from the coefficients of its normal equation". (Russian)  Mat. Zametki  73  (2003),  no. 3, 453--456;  translation in  Math. Notes  73  (2003),  no. 3-4, 419--423

\bibitem{kruch} Mubarakzyanov, G. M. "On solvable Lie algebras", (Russian) Izv. Vis. Ucheb. Zav. Matem. (1963)  no 1 (32) P. 114-123.

\bibitem{nester} M.Nesterenko, R.Popovych, "Contractions of low-dimensional Lie algebras".  J. Math. Phys.  47  (2006),  no. 12, 123515, 45 pp.


\bibitem{poincare} Poincare\,H. "Les fonctions analytiques de deux variables et la
representation conforme",
Rend.\,Circ.\,Mat.\,Palermo.\,1907.\,23.\,P.185-220.

\bibitem{nester1} Popovych, R.; Boyko, V.; Nesterenko, M.; Lutfullin, M. "Realizations of real low-dimensional Lie algebras".  J. Phys. A  36  (2003),  no. 26, 7337--7360.

\bibitem{shirokov} Shirokov, P. A.; Shirokov, A. P. {Affinnaya differentsialnaya geometriya.} (Russian) [Affine differential geometry] Gosudarstv. Izdat. Fiz.-Mat. Lit., Moscow 1959 319 pp.

\bibitem{tomassini} Tomassini, G., Tracce delle funzioni olomorfe sulle sottovariet\`a analitiche reali d'una variet\`a complessa, {\it Ann. Scuola Norm. Sup. Pisa} (3) 20(1966), 31--43.

\bibitem{zaitsev} Zaitsev, D., "On different notions of homogeneity for CR-manifolds",  Asian J. Math.  11  (2007),  no. 2, 331--340.



\end{thebibliography}
\end{document}